\documentclass[reqno,12pt]{amsart}
\newtheorem{theorem}{Theorem}[section]
\newtheorem{lemma}[theorem]{Lemma}
\newtheorem{corollary}[theorem]{Corollary}

\theoremstyle{definition}
\newtheorem{assumption}[theorem]{Assumption}
\newtheorem{definition}[theorem]{Definition}
\newtheorem{exercise*}[theorem]{Exercise*}

\theoremstyle{remark}
\newtheorem{example}[theorem]{Example}
\newtheorem{remark}[theorem]{Remark}

\makeatletter
\def\dashintindex{\operatorname%
{-\kern-.7em\DOTSI\intop\ilimits@}}%
\def\dashint{\operatorname%
{\,\,\text{\bf--}\kern-.98em\DOTSI\intop\ilimits@\!\!}}
\def\dashnorm{\,\,\text{\bf--}\kern-.5em\|}
\def\ninf{\qopname\relax\@empty{inf\phantom{p}\!\!\!}}
\makeatother

\newcommand\cbrk{\text{$]$\kern-.15em$]$}} 
\newcommand\opar{\text{\raise.2ex\hbox{${\scriptstyle | }$}\kern-.34em$($} }

\newcommand\frc{\mathfrak{c}}
\newcommand{\cedc}{\hbox{\rm \c{c}}}

\newcommand\bB{\mathbb{B}}

\newcommand\bR{\mathbb{R}}

\def\sff{{\sf f}}

\newcommand\cB{\mathcal{B}}

\newcommand\cF{\mathcal{F}}

\newcommand\cP{\mathcal{P}}

\newcommand{\sign}{\text{\rm\,sign}\,}

 \newcommand{\mysection}[1]{\section{#1}
 \setcounter{equation}{0}}

\begin{document}

\title[$L_{p}$-estimates for solutions of SPDEs]
{Estimates in $L_{p}$ for solutions of SPDEs with coefficients in Morrey classes}
\author{N.V. Krylov}
\email{nkrylov@umn.edu}
\address{127 Vincent Hall, University of Minnesota, Minneapolis, MN, 55455}



\begin{abstract}
For solutions of a certain class of SPDEs
in divergence form we present some estimates of their
$L_{p}$-norms and the $L_{p}$-norms of their
first-order derivatives. The main novelty
is that the low-order coefficients
are supposed to belong to certain Morrey
classes instead of $L_{p}$-spaces. Our results are new
even if there are no stochastic terms in the equation.
\end{abstract}

\maketitle

\mysection{Introduction}

In this paper we come back to an old problem
of estimating the $L_{p}$-norm in $x$ of
solutions $u_{t}(x)$ of It\^o stochastic 
partial differential equations of second order.
The first such estimates appeared in \cite{KR_77} in 1977 and were achieved by using It\^o's formula for $\|u_{t}\|^{p}_{L_{p}}$ and integrating by parts. This method of obtaining such estimates was generalized, extended, and applied in many papers of which we cite, probably, the most recent 
\cite{DG_15}, \cite{GGK_15}, \cite{Qi_20}, and 
\cite{TWT_20} containing vast lists of references.
The closest to the present article are
some computations in \cite{Kr_304} which
is about strong solutions of It\^o stochastic equations rather than about SPDEs. 
There, as in the present article, the main emphasis is on the drift coefficients
from a Morrey class.
In this connection it is also worth mentioning
\cite{KM_21} and some   references therein  where methods different from integration by parts are
used to treat strong solutions of It\^o stochastic equations with singular drift. We apply integration by parts to SPDEs which may be just nonrandom usual parabolic equations. Even in this case,
as far as the author can tell, the presented results are new although similar results
for the {\em elliptic\/} equations can be found in
\cite{MT_15}.

Let $d,d_{1}$ be integers, $d\geq 3$, $\bR^{d}$ be a $d-$dimensional Euclidean space of points
$x=(x^{1},...,x^{d})$, and
$$
\bR^{d+1}=\{(t,x):t\in\bR,x\in\bR^{d}\}.
$$
Let $\cB(\bR^{d})$ be the Borel $\sigma$-field on $\bR^{d}$, $L_{p}=L_{p}(\bR^{d})$.

Let $(\Omega,\cF,P)$ be a complete probability space 
 equipped with a filtration 
$(\cF_t)_{t\geq0}$   of complete $\sigma$-fields
$\cF_{t}\subset\cF$. We suppose that
on our probability space we are given a $d_{1}$-dimensional   Wiener process $w_{t}$, which is a Wiener process relative to the filtration 
$(\cF_t)_{t\geq0}$. By $\cP$ we denote
the predictable $\sigma$-field on $\Omega\times(0,\infty)$ (we follow the terminology in \cite{IW}). Fix a bounded stopping time $\tau$.
 
Assume that on $\Omega\times(0,\infty)
\times\bR^{d}$ we are given 
$\cP\otimes\cB(\bR^{d})$-measurable $\bR^{d}$-valued functions
  $\sigma^{k}_{t}(x)=(\sigma^{ik}_{t}(x))$,
$k=1,...,d_{1}$, $\beta_{t}(x)=(\beta^{i}_{t}(x))$,
$ b _{t}(x)=( b ^{i}_{t}(x))$,
$ \sff_{t}(x)=( \sff^{i}_{t}(x))$,
  $d\times d$-matrix valued
functions $a_{t}(x)= (a_{t}^{ij}(x))$, 
$\bR^{d_{1}}$-valued functions
$\nu _{t}(x)=(\nu^{k}_{t}(x))$,
$g _{t}(x)=(g^{k}_{t}(x))$,
and real-valued
functions $c_{t}(x)$, $f_{t}(x)$. Introduce
$$
L_{t}u_{t}=D_{i}(a^{ij}_{t}D_{j}u_{t}+\beta^{i}_{t}u_{t})+b^{i}_{t}
D_{i}u_{t}+c_{t}u_{t},
\quad M^{k}_{t}u_{t}=\sigma^{ik}_{t}
D_{i}u_{t}+\nu^{k}_{t}u_{t},
$$
where and everywhere below the common summation 
convention is used and $D_{i}=\partial/\partial x^{i}$.

 Our goal is to investigate the equation
\begin{equation}                                                       
                                 \label{12.29.1}
du_{t}=(L_{t}u_{t}+D_{i}\sff^{i}_{t}+f_{t})\,dt+(M^{k }_{t}u_{t}+g^{k }_{t})\,dw^{k }_{t},\quad
t\leq \tau,
\end{equation}
with initial condition
\begin{equation}                                                       
                                 \label{12.29.2}
u_{0}=\Phi.
\end{equation}

\begin{remark}
                   \label{remark 1.5.1}

In the literature very popular conditions
on $ b$ (and $\beta$) is that
$b\in L_{p_{0},q_{0}}(T)$, that is, for   $T\in(0,\infty)$,
\begin{equation}
                            \label{1.5.1}
\Big(\int_{0}^{T}\Big(\int_{\bR^{d}}|b_{t}(x)|^{p_{0}}
\,dx\Big)^{q_{0}/p_{0}}\,dt\Big)^{1/q_{0}}<\infty
\end{equation}
with $p_{0},q_{0}\in[2,\infty]$ satisfying
$$
\frac{d}{p_{0}}+\frac{2}{q_{0}}=1.
$$ 

Observe that, if $p_{0}>d$ and we take an
arbitrary constant $\hat N$
and introduce
$$
\lambda(t)=\hat N\Big(\int_{\bR^{d}}
|b_{t}(x)|^{p_{0}}\,dx\Big)^{1/(p_{0}-d )},
$$
then for (the notation will be explained later)
$$
b^{M}_{t}( x)=b_{t} ( x)I_{|b_{t} (x)|\geq \lambda(t)}
$$
 we have
$$
\dashint_{B_{\rho}}|b^{M}_{t}( x)|^{d}\,dx
\leq \lambda^{d-p_{0}}(t)
\dashint_{B_{\rho}}|b_{t}  ( x)|^{p_{0}}\,dx
\leq N(d)\hat N^{d-p_{0}}\rho^{-d},
$$
where $B_{\rho}$ is a ball of radius $\rho
\in(0,\infty)$,  
$$
\dashint_{B_{\rho}}...=\frac{1}{|B_{\rho}|}\int_{B_{\rho}}...,\quad |B_{\rho}|=\text{\rm Vol}\,
B_{\rho}.
$$

Here $N(d)\hat N^{d-p_{0}}$ can be made
arbitrarily small if we choose $\hat N$
large enough. In addition, for 
$b^{B}_{t}=b_{t}-b^{M}_{t}$ we have
 $|b^{B}_{t}|\leq \lambda(t)$ and
$$
\|b^{B}\|^{2}_{L_{\infty,2}(T)}\leq
\int_{0}^{T}\lambda^{2}(t)\,dt
=\hat N^{2}\int_{0}^{T}\Big(\int_{\bR^{d}}|b_{t}(x)|^{p_{0}}
\,dx\Big)^{q_{0}/p_{0}}\,dt<\infty.
$$
This shows that the assumptions we are going to
impose on $b$ are weaker then \eqref{1.5.1}
if $p_{0}>d$.

If $p_{0}=d$ (and $q_{0}=\infty$) our assumptions are weaker if  \eqref{1.5.1} is combined
with the requirement that the family
$|b_{t}|^{d},t\in(0,\infty)$, be uniformly integrable (say, case of $b$ independent of $t$).
\end{remark}
 
Recall that $W^{n}_{p}$ are Sobolev spaces
of functions $u$ on $\bR^{d}$ such that $u$
and all its (generalized) derivatives
of order $\leq n$
belong to $L_{p}$. The norm in $W^{n}_{p}$
is introduced in a natural way.
 
\begin{definition}                                                               \label{definition solution}
A $W_p^1$-valued function $u$,
defined on the stochastic interval $\opar0,\tau\cbrk$, is called a solution 
of 
\eqref{12.29.1}-\eqref{12.29.2} 
  if $u$ is predictable on $\opar0,\tau\cbrk$, 
$$
\int_0^{\tau}\|u_t\|_{W^1_p}^p\,dt<\infty , 
$$
(a.s.) and for each $\varphi\in C_0^{\infty}(\bR^d)$ for almost all $\omega\in\Omega$
$$
(u_t,\varphi)= (\Phi,\varphi)+\int_0^t(\sigma^{ik}_sD_iu_s+\nu^{k}_{s}u_s+g^k_s,\varphi)\,dw^k_{s}
$$
\begin{equation}
                            \label{1.11.2}
+\int_0^t\{( b^{i}_sD_iu_s+c_su_s+f_s,\varphi)-(a^{ij}_sD_ju_s+\beta^{i}_{s}u_{s}+\sff^{i}_{s},D_i\varphi)
\}\,ds
\end{equation}
for all $t\in[0,\tau(\omega)]$,   where   $(u,v)$ denotes the  integral over  $  \bR^d $   of $uv$.
\end{definition}

Our standing assumption is that equation \eqref{12.29.1} is uniformly nondegenerate.
Set $\alpha_{t}(x)=(\alpha_{t}^{ij}(x))$,
where $\alpha_{t}^{ij} =\sigma^{ik}_{t}\sigma^{jk}_{t}$. Note that we do not suppose that $a_{t}$
is symmetric.
\begin{assumption}
                      \label{assumption 12.29.1}
There is a constant $\delta>0$ such that
for all values of arguments and $\lambda\in\bR^{d}$
\begin{equation}
                               \label{12.29.3}
|a_{t}|\leq\delta^{-1},\quad
(2a^{ij}_{t}-\alpha_{t}^{ij})\lambda^{i}\lambda^{j}\geq \delta|\lambda|^{2}.
\end{equation}
\end{assumption}

Another standing assumption concerns free
terms. We fix $p\in[2,d)$ ($d\geq3$) and recall
that
$L_{p}=L_{p}(\bR^{d})$.
\begin{assumption}
                      \label{assumption 12.29.2}
We have
\begin{equation}
                               \label{12.29.4}
E\int_{0}^{\tau}(\|f_{t}\|^{p}_{L_{p}}+
\|\sff_{t}\|^{p}_{L_{p}}+\|g_{t}\|^{p}_{L_{p}})
\,dt<\infty.
\end{equation}
\end{assumption}

The assumptions on the coefficients of $L$ and $M^{k}$ is more delicate.
Fix $r\in (p, d] $, $R_{0}\in(0,\infty)$,
and let $\bB_{\rho}$ denote the set of balls in $\bR^{d}$ with radius $\rho$.

\begin{definition}
                  \label{definition 12.16.1}
We call a real- or vector- or else tensor-valued
function $f_{t}(x)$ given on $\Omega\times(0,\infty)
\times\bR^{d} $ {\em admissible\/} if it  is
represented as $f_{t}= f^{M}_{t}+f^{B}_{t}$ (``Morrey part'' of
$f$ plus its ``bounded part'')
with  
$\cP\otimes\cB(\bR^{d})$-measurable
$f^{M}_{t}$, $f^{B}_{t}$   such that there exists a {\em constant\/} $\hat f<\infty$ for which
$$
\Big(\dashint_{B } |f^{M}_{t}|^{ r}\,dx\Big)^{1/r}\leq \hat f  \rho^{-1},
$$
whenever $B\in\bB_{\rho}$ and $\rho\leq R_{0}$,
and there exists a {\em predictable\/} $\bar f _{t}\geq
\sup_{x\in\bR^{d}}|f^{B}_{t}( x)|$
such that the integral
$$
\int_{0}^{\tau}|\bar f _{t} | ^{2}\,dt 
$$
is a bounded function on $\Omega$.
\end{definition}

\begin{assumption}
                  \label{assumption 12.19.1}
$ \beta,b,\cedc:= c_{+}^{1/2}, \nu $ are admissible.

\end{assumption}

Here is our first main result about an a priori estimate in which for convenience
of   applications we introduce
  a measurable, $\cF_{t}$-adapted process
$\mu_{t}\geq 0$ such that the integral
$$
\int_{0}^{\tau} \mu_{t} \,dt 
$$
is a bounded function on $\Omega$.

\begin{theorem}
                   \label{theorem 12.19.1}
Under the above assumptions 
suppose that a solution $u_{t}$ of \eqref{12.29.1}-\eqref{12.29.2} exists.
  Then with probability one $u_{t\wedge\tau}$
is continuous as an $L_{p}$-valued function and
there exists a constant $N_{0}=N_{0}(p,\delta,d,r)$ such that if
\begin{equation}
                           \label{12.19.20}
N_{0}(  
\hat b+\hat\beta+\widehat{\cedc}+\hat \nu )\leq 1,
\end{equation}
then for  any stopping time $\gamma\leq\tau$ 
$$
Ee^{-\phi_{\gamma}}\int_{\bR^{d}}|u_{t}|^{p}\,dx
+(\delta/16)E\int_{0}^{\gamma}e^{-\phi_{s}}
\int_{\bR^{d}}|DU_{s}|^{2}\,dxds
$$
\begin{equation}
                           \label{12.19.2}
\leq
E\int_{\bR^{d}}|u_{0}|^{p}\,dx +NE\int_{0}^{\gamma}e^{-\phi_{s}}
\int_{\bR^{d}}(|\sff_{t}|^{p}+|f_{t}|^{p}+|g_{t}|^{p}
-\mu_{t}|u_{t}|^{p})\,dxds,
\end{equation}
$$
E\sup_{t\leq\tau}e^{-\Phi_{t}}\int_{\bR^{d}}|u_{t}|^{p}
\,dx\leq N\int_{\bR^{d}}|u_{0}|^{p}\,dx
$$
\begin{equation}
                         \label{1.10.1}
+
NE\int_{0}^{\tau}e^{-\Phi_{s}}
\int_{\bR^{d}}( |\sff_{s}|^{p}+|f_{s}|^{p}+|g_{s}|^{p} )\,dxds,
\end{equation}
where $U_{t}=|u_{t}|^{p/2}$,
$$
\Phi_{t}=(B+1)\int_{0}^{t} \lambda_{s}\,ds,\quad\phi_{t}=\int_{0}^{t}(B\lambda_{s}+\mu_{s})\,ds,\quad \lambda_{t}= 1+\bar \beta^{2}_{t}+\bar b_{t}^{2}+\bar \cedc_{t}^{2}
+\bar \nu_{t}^{2} ,
$$
and the constants $N, B$ depend only on $p,r,\delta,d, R_{0}$. 

\end{theorem}
We prove this theorem in Section \ref{section 1.14.2}.

\begin{example}[Example 1.3.1 \cite{LSU_67}]
                       \label{example 1.5.1}
Condition \eqref{12.19.20} requires, in particular, $\hat b$ to be sufficiently small.
It turns out that this smallness assumption
is essential. For instance, $u_{t}(x)=\exp(-|x|^{2}/(4t))$ satisfies
\begin{equation}
                             \label{1.9.10}
\frac{\partial}{\partial t}u_{t}(x)
=\Delta u_{t}(x)-\frac{d}{|x|^{2}}x^{i}D_{i}u_{t},
\end{equation}
where $|b_{t}(x)|=d/|x|$ is admissible
with any $r<d$. In this situation \eqref{12.19.2} fails. However, if we take
$\varepsilon b$ in place of $b$ our  equation
will satisfy \eqref{12.19.20} and the estimate
becomes available if $\varepsilon$
is sufficiently small.

Similarly, one cannot allow $\hat\beta$ and $\hat\cedc$ to be large. Indeed, equation
\eqref{1.9.10} can be rewritten as
$$\frac{\partial}{\partial t}u_{t}(x)
=\Delta u_{t}(x)-D_{i}\Big(\frac{d}{|x|^{2}}x^{i} u_{t}\Big)+\frac{d(d-2))}{|x|^{2}}u_{t}.
$$

Note that $|b|=\varepsilon/|x|$ is not allowed
in \cite{LSU_67} and \cite{MT_15}, so it seems
that Theorem \ref{theorem 12.19.1} is new
even if there are no stochastic terms.

\end{example}

Our next result is about estimating
$Du_{t}$. For that we need to impose a stronger assumptions than Assumption \ref{assumption 12.29.2} and \ref{assumption 12.19.1}.
By $Df(x)$ we mean the collection of the
first-order Sobolev derivatives of
(vector- or   tensor-valued) $f(x)$ and by
$|Df(x)|$ we mean any fixed norm in the space
where $Df(x)$ takes its values.
Similar meaning will be given later to $D^{2}f$
and $|D^{2}f|$.

\begin{assumption}
                     \label{assumption 1.5.2}
We have $\beta=\sff\equiv 0$ and
$$
E\int_{0}^{\tau} 
\int_{\bR^{d}}(|f_{t}|^{p}+|g_{t}|^{p}+|Df_{t}|^{p}+|Dg_{t}|^{p})\,dxds<\infty.
$$

\end{assumption}

Introduce $\frc_t=|Dc_{t}|^{1/2}$.
\begin{assumption}
                  \label{assumption 12.20.1}
The functions
$Da, D\sigma, b,\cedc ,\frc, \nu,D\nu $ are admissible.
\end{assumption}

\begin{theorem}
                   \label{theorem 12.24.1}
Under the above assumptions
  suppose that a solution $u_{t}$ of \eqref{12.29.1}-\eqref{12.29.2} exists
and
\begin{equation}
                           \label{1.10.6}
E\int_{0}^{\tau}\|u_{t}\|_{W^{2}_{p}}^{p}\,dt<\infty.
\end{equation}
 Then (a.s.) $Du_{t\wedge \tau}$ is a continuous $L_{p}$-valued function and
there exists a constant $N_{0}=N_{0}(p,\delta,d,r)$ such that, if
\begin{equation}
                           \label{12.24.10}
N_{0}(  \widehat{Da}+\widehat{D\sigma}+
\hat b+\widehat{\cedc}+\hat\frc+\hat \nu +\widehat{D\nu} )\leq 1,
\end{equation}
then   
\begin{equation}
                           \label{12.24.9}
G+H\leq N K,
\end{equation}
where  
$$
G=E\sup_{t\leq\tau}e^{-\Psi_{t}}\int_{\bR^{d}}|Du_{t}|^{p}\,dx,
$$
$$
H= E\int_{0}^{\tau}e^{-\Psi_{s}}
\int_{\bR^{d}}\big(|Du_{t}|^{p-2}
|D^{2}u_{t}|^{2}+\Lambda_{t}|Du_{t}|^{p}\big)\,dxds,
$$ 
$$
K= E\int_{\bR^{d}} ( |u_{0}|^{p}+|Du_{0}|^{p})\,dx
$$
$$
+E\int_{0}^{\tau}e^{-\Psi_{s}}
\int_{\bR^{d}}(|f_{t}|^{p}+|g_{t}|^{p}+|Df_{t}|^{p}+|Dg_{t}|^{p})\,dxds ,
$$
$$
\Psi_{t}=C \int_{0}^{t}\Lambda_{s}\,ds,\quad
\Lambda_{t} = 1+\overline{Da}^{2}_{t}+
\overline{D\sigma}_{t}^{2}+
\bar b_{t}^{2}+\bar \cedc_{t}^{2}+\bar \frc_t^2
+\bar \nu_{t}^{2}+\overline{D\nu_{t}}^{2} 
$$
and the constants $N,C$ depend only on $p,r,\delta,d, R_{0}$. 

\end{theorem}

We prove this theorem in Section \ref{section 1.14.1}.

\begin{remark}
Somewhat similar result is presented in
Theorem 2.1 of \cite{GGK_15}. There are, however, several distinctions. Let alone systems, in 
\cite{GGK_15} the equations can degenerate
and because of that the assumptions on
the smoothness of $a$ and 
$\sigma$ are stronger. Also the coefficients of lower order terms like $b,c$ 
are assumed to be bounded. We heavily use
the nondegeneracy of our equation.

For nondegenerate equations \cite{Kr_09}
provides a similar result with a weaker assumption on $a$ but with bounded lower order coefficients  and with an  assumption on $\sigma$ which, basically, requires its continuity in $x$. Under our assumptions $\sigma$ can have rather wild discontinuities.
Incidentally, one of the aims of this article
is to develop some tools on the way of proving that under our assumptions on $\sigma$ and $b$ the stochastic equation
$dx_{t}=\sigma^{k}_{t}(x_{t})\,dw^{k}_{t}+
b_{t}(x_{t})\,dt$ with deterministic $\sigma$
and $b$ and any nonrandom initial condition has a  
strong solution provided that the matrix
$(\sigma^{ik}_{t}\sigma^{jk}_{t})$ is uniformly
nondegenerate. The way we mean is similar to the one 
in \cite{Kr_304}.

\end{remark}

\mysection{Auxiliary results}

For   $r\in(1,d]$ and
  two real-valued functions
$b$ and $\cedc$ on $\bR^{d}$ set
$$
\hat b:=\|b\|_{E_{r,1}} ,\quad
\hat \cedc:=\|\cedc\|_{E_{r ,1}} ,
$$
where, for any function $f$,
$$
\|f\|_{E_{r,1}}:=\sup_{\substack{\rho>0\\B\in\bB_{\rho}}}\rho \dashnorm f\|_{L_{r}(B)},
\quad \dashnorm f\|_{L_{r}(B)}^{r}
=\dashint_{B}|f|^{r}\,dx.
$$
For functions $u,v$ on $\bR^{d}$
we write $u\prec v$ if the integral of $u$
over $\bR^{d}$ is less than or equal to that
of $v$. If the integrals are equal, we write
$u\sim v$. Also set $p'=p/(p-1)$.
 
\begin{lemma}
                      \label{lemma 1.6.1}
Let $u,v \in C^{\infty}_{0}$, then
\begin{equation}
                          \label{1.9.1}
1<p<r \Longrightarrow|bu|^{p}\prec N\hat b^{p}
|Du|^{p},
\end{equation}
\begin{equation}
                          \label{1.6.2}
  p,p'\in (1,r) \Longrightarrow
\int_{\bR^{d}}|v\cedc^{2} u|\,dx
\leq N\hat \cedc^{2}\|Dv\|_{L_{p'}}
\|Du\|_{L_{p}},
\end{equation}
\begin{equation}
                          \label{1.6.3}
1< p'<r  \Longrightarrow
\int_{\bR^{d}}|vbDu|\,dx
\leq N\hat b\|Dv\|_{L_{p'}}
\|Du\|_{L_{p}},
\end{equation}
where the constants $N$ depend only on $d,r,p$.
\end{lemma}

Proof. Estimate \eqref{1.9.1}
is proved in \cite{CF_90}. To prove \eqref{1.6.2},
it suffices to observe that
$$
\int_{\bR^{d}}|v\cedc^{2}u|\,dx
\leq\|\cedc v\|_{L_{p'}}\|\cedc u\|_{L_{p}},
$$
and then use \eqref{1.9.1}.
Similarly \eqref{1.6.3} is proved on the basis
of
$$
\int_{\bR^{d}}|vbDu|\,dx\leq\|bv\|_{L_{p'}}
\|Du\|_{L_{p}}.
$$
 The lemma is proved.

\begin{remark}
                       \label{remark 1.11.1}
Estimate \eqref{1.6.3} shows, in particular,
that $g:=b^{i}D_{i}u $, with $\widehat {|b|}<\infty$,    can be represented as $D_{i}\sff^{i}+f$, where the $L_{p}$-norms of $\sff^{i},f$
are controlled by that of $Du$.
Indeed, for $\Lambda=(1-\Delta)^{1/2}$
the estimate implies that $h:=\Lambda^{-1}g\in L_{p}$ and $g=\Lambda h=(1-\Delta)\Lambda^{-1}h=D_{i}\sff^{i}+f$, where $f=\Lambda^{-1}h$,
$\sff^{i}=-D_{i}\Lambda^{-1}h$.

\end{remark}
 
\begin{lemma}
                          \label{lemma 12.17.1}
If $f$ is admissible, then for any $t$ and $u
\in C^{\infty}_{0}$ we have
\begin{equation}
                               \label{12.17.1}
|f_{t}|^{2}|u|^{2}\prec N(d,r)\hat f^{2}|Du|^{2}
+(N(d,r)R_{0}^{-2}\hat f^{2}+2  \bar f^{2}_{t} )|u|^{2}.
\end{equation}
\end{lemma}

Proof. We have $|f_{t}|^{2}\leq 2|f^{M}_{t}|^{2}+2|f^{B}_{t}|^{2}$ and this shows how
$2 |\bar f^{B}_{t}|^{2} |u|^{2}$ appeared in 
\eqref{12.17.1}. By a version of the
Chiarenza--Frasca theorem in \cite{CF_90}
(see, for instance, Lemma 3.5 in
\cite{Kr_304})
$$
|f^{M}_{t}|^{2}|u_{t}|^{2} 
\prec N \hat f^{2}|Du_{t}|^{2}
+N R_{0}^{-2 }\hat f^{2} 
 | u_{t}|^{2}.
$$
This proves the lemma.
\begin{corollary}
                     \label{corollary 12.17.1}
$$
|f_{t}|\,|Du_{t}|\,|u_{t}|\leq \varepsilon|Du_{t}|^{2}
+\varepsilon^{-1}|f_{t}|^{2}|u_{t}|^{2} 
$$
$$
\prec (\varepsilon+N \varepsilon^{-1}\hat f^{2})|Du_{t}|^{2}+N\varepsilon^{-1}(\bar f^{2}_{t}+ R_{0}^{-2 }\hat f^{2} )
 | u_{t}|^{2}.
$$
\end{corollary}

\mysection{Proof of Theorem \ref{theorem 12.19.1}}

                     \label{section 1.14.2}

  In light of Lemma \ref{lemma 1.6.1}
and Remark \ref{remark 1.11.1},
Theorem 2.1 of \cite{Kr_10} is applicable,
which in particular implies the $L_{p}$-continuity of $u_{t\wedge\tau}$. By this theorem
also, after introducing   $\theta_{t}
=\sign u_{t}$ we obtain for $t\leq\tau$ that
$$
(1/p)\int_{\bR^{d}}|u_{t}|^{p}\,dx=
(1/p)\int_{\bR^{d}}|u_{0}|^{p}\,dx
+\int_{0}^{t}\int_{\bR^{d}}J_{s}\,ds+ m_{t}
$$
$$
+\int_{0}^{t}\int_{\bR^{d}}\theta_{s}|u_{s}|^{p-1} 
(L_{s}u_{s}+D_{i}\sff^{i}_{s}+f_{s})\,dx\,ds,
$$
where
$$
J_{t}:=((p-1)/2)|u_{s}|^{p-2}\sum_{k}\big(M^{k }_{s}u_{t}+g^{k }_{s})\big)^{2}\,dx
$$
$$
\leq ((p-1)/2)|u_{t}|^{p-2}\sigma^{ik}_{t}
D_{i}u_{t}\sigma^{jk}_{t}
D_{j}u_{t}+(\delta/(2p))|DU_{t}|^{2}
$$
\begin{equation}
                              \label{1.10.3}
+N|u_{t}|^{p}|\nu_{t}|^{2}+N|u_{t}|^{p-2}
|g_{t}|^{2},
\end{equation}
$$
m_{t}:=p\int_{0}^{t}\int_{\bR^{d}}\theta_{s}|u_{s}|^{p-1}\big(M^{k }_{s}u_{s}+g^{k }_{s})\big)\,dx\,dw^{k}_{s}.
$$
 Here and below  constants $N $ denote
constants depending only on $d,r,p,\delta$.
 
From the start we may assume that
$\hat b,\hat\beta,\hat\cedc,\hat\nu \leq1$
(in order to be able to drop terms like $\hat f^{2}$ in the second term on the right in \eqref{12.17.1}). 
Note that
$$
I_{1}:=\theta_{t}|u_{t}|^{p-1}D_{i}(a^{ij}_{t}D_{j}
u_{t})+((p-1)/2)|u_{t}|^{p-2}\sigma^{ik}_{t}
D_{i}u_{t}\sigma^{jk}_{t}
D_{j}u_{t}
$$
$$
+(\delta/(2p))|DU_{t}|^{2}\sim-\frac{2(p-1)}{p^{2}}(2a-\alpha)^{ij}_{t}
D_{i}U_{t}D_{j}U_{t}
$$
\begin{equation}
                              \label{12.19.4}
+(\delta/(2p))|DU_{t}|^{2}\leq -\hat \delta|DU_{t}|^{2},
\end{equation}
where $\hat \delta= \delta/(2p) $.

Next, by Corollary \ref{corollary 12.17.1}
$$
I_{2}:=\theta_{t}|u_{t}|^{p-1}D_{i}(\beta^{i}_{t}u_{t}) 
+\theta_{t}|u_{t}|^{p-1}b^{i}_{t}D_{i}u_{t}
$$
$$
\sim-(p-1)\theta_{t}|u_{t}|^{p-1}\beta^{i}_{t}D_{i}u_{t}+\theta_{t}|u_{t}|^{p-1}b^{i}_{t}D_{i}u_{t}
$$
$$
\leq (NU_{t}(|\beta_{t}|+|b_{t}|)\big)|DU_{t}|
$$
$$
\leq (\hat \delta/4+N_{1}(\hat \beta^{2}+\hat b^{2}))|DU_{t}|^{2}+N (1+\bar \beta^{2}_{t}+\bar b^{2}_{t}  )
 | U_{t}|^{2}.
$$
 
We subject $\hat \beta$ and $\hat b$ to
\begin{equation}
                              \label{12.19.5}
N_{1}(\hat \beta^{2}+\hat b^{2})\leq \hat\delta/4
\end{equation}
and then get that
\begin{equation}
                              \label{12.19.6}
I_{2}
\leq (\hat \delta/2 )|DU_{t}|^{2}+N (1+\bar \beta^{2}_{t}+\bar b^{2}_{t}  )
 | u_{t}|^{p}.
\end{equation}

By Lemma \ref{lemma 12.17.1}
$$
I_{3}:=\theta_{t}|u_{t}|^{p-1}c_{t}u_{t}+N|u_{t}|^{p}
|\nu_{t}|^{2}\leq
\cedc_{t}^{2}U_{t}^{2}+N|\nu_{t}|^{2}U_{t}^{2}
$$
$$
\prec N_{2}(\hat \cedc^{2}+\hat \nu^{2})|DU_{t}|^{2}+
N(1+\bar\cedc_{t}^{2}+\bar \nu^{2}_{t})|U_{t}|^{2}.
$$
We subject $\hat \cedc$ and $\bar\nu$ to
\begin{equation}
                              \label{12.19.7}
N_{2}(\hat \cedc^{2}+\hat \nu^{2})\leq \hat\delta/4
\end{equation}
and then get that
\begin{equation}
                              \label{12.19.8}
I_{3}
\leq (\hat \delta/4 )|DU_{t}|^{2}+N (1+\bar \cedc^{2}_{t}+\bar \nu^{2}_{t} )
 | u_{t}|^{p}.
\end{equation}
Next,
 $$
I_{4}:=|u_{t}|^{p-1}\theta_{t}D_{i}\sff^{i}_{t}
\sim N \sff^{i}_{t}|u_{t}|^{p/2-1}D_{i}U_{t}
\leq (\hat\delta/16)|DU_{t}|^{2}+N
|\sff_{t}|^{2}|u_{t}|^{p-2} 
$$
\begin{equation}
                              \label{1.9.2}
\leq (\hat \delta/8)|DU_{t}|^{2}+N|U_{t}|^{2}
+N |\sff_{t}|^{p}.
\end{equation}

In estimating the last remaining terms
we use elementary inequalities
\begin{equation}
                            \label{1.13.1}
 |u_{t}|^{p-1}f_{t}+N|u_{t}|^{p-2}|g_{t}|^{2}\leq N|u_{t}|^{p}+N(|f_{t}|^{p}+|g_{t}|^{p}).
\end{equation}
Upon combining this with  
\eqref{12.19.4}, \eqref{12.19.6},
\eqref{12.19.8}, and \eqref{1.9.2}, we get
that, if condition \eqref{12.19.20}
is satisfied for an appropriate $N_{0}$, then
$$
(1/p)d\int_{\bR^{d}}|u_{t}|^{p}\,dx\leq 
\int_{\bR^{d}}\Big(-(\hat\delta/8)|DU_{t}|^{2}+N(|\sff_{t}|^{p}+|f_{t}|^{p}+|g_{t}|^{p})\Big)\,dx\,dt
$$
\begin{equation}
                         \label{1.10.2}
+(B/p)\int_{\bR^{d}}\lambda_{t}|u_{t}|^{p}\,dx\,dt
+dm_{t}.
\end{equation}
It follows that for $t\leq\tau$
$$
e^{-\phi_{t}}\int_{\bR^{d}}|u_{t}|^{p}\,dx
+(\hat\delta/8)\int_{0}^{t}e^{-\phi_{s}}
\int_{\bR^{d}}|DU_{s}|^{2}\,dxds
$$
$$
\leq \int_{\bR^{d}}|u_{0}|^{p}\,dx+N\int_{0}^{t}e^{-\phi_{s}}
\int_{\bR^{d}}( |\sff_{s}|^{p}+|f_{s}|^{p}+|g_{s}|^{p}-\mu_{t}|u_{t}|^{p})\,dxds
$$
$$
+\int_{0}^{t}e^{-\phi_{s}}\,dm_{s}.
$$
 This yields \eqref{12.19.2} because,
as we know from \cite{Kr_10},
$m_{t}$ is a martingale.
To prove \eqref{1.10.1} we follow
E. Pardu. Estimate \eqref{1.10.2} and the Davis's
inequality imply that
$$
E\sup_{t\leq\tau}e^{-\Phi_{t}}\int_{\bR^{d}}|u_{t}|^{p}
\,dx\leq E\int_{\bR^{d}}|u_{0}|^{p}\,dx
$$
$$
+
NE\int_{0}^{\tau}e^{-\Phi_{s}}
\int_{\bR^{d}}( |\sff_{s}|^{p}+|f_{s}|^{p}+|g_{s}|^{p} )\,dxds
$$  
\begin{equation}
                             \label{1.10.5}
+NE\Big(\int_{0}^{\tau}e^{-2\Phi_{s}}\sum_{k}\Big(
\int_{\bR^{d}}
|u_{s}|^{p-1}|M^{k} _{s}u_{s}+g^{k} _{s}|\,dx\Big)^{2}\,ds\Big)^{1/2}.
\end{equation}
Here $|u_{s}|^{p-1}=|u_{s}|^{p/2 }|u_{s}|^{p/2-1}$ and the last term is dominated by
$$
I:=NE\Big(\int_{0}^{\tau}e^{-2\Phi_{s}}\|u_{s}\|_{L_{p}}^{p} 
\int_{\bR^{d}} 
J_{s}\,dx \,ds\Big)^{1/2}
$$
$$
 \leq NE\Big(\sup_{s\leq\tau}e^{-\Phi_{s}}\|u_{s}\|_{L_{p}}^{p }\Big)^{1/2}\Big( \int_{0}^{\tau}e^{-\Phi_{s}}\int_{\bR^{d}} 
J_{s}\,dx \,ds\Big)^{1/2}
$$
$$
\leq N\Big(E\sup_{s\leq\tau}e^{-\Phi_{s}}\|u_{s}\|_{L_{p}}^{p }\Big)^{1/2}
\Big( E\int_{0}^{\tau}e^{-\Phi_{s}}\int_{\bR^{d}} 
J_{s}\,dx \,ds\Big)^{1/2}
$$
It is seen from \eqref{1.10.3} and the estimates of $I_{3}$ and \eqref{1.13.1} that
$$
J_{t}\leq N|DU_{t}|^{2}+
N|u_{t}|^{p}|\nu_{t}|^{2}+N|u_{t}|^{p-2}|g_{t}|^{2}
$$
$$
\prec N|DU_{t}|^{2}+N(1+\bar \nu^{2}_{t})
|u_{t}|^{p}+N|g_{t}|^{p}.
$$
The last expression is under control
due to \eqref{12.19.2} with
$\mu_{t}=\lambda_{t}\geq 1+\bar \nu^{2}_{t}$. Hence coming back to \eqref{1.10.5} 
we immediately obtain \eqref{1.10.1}.
This proves the theorem.

\begin{remark}
                        \label{remark 12.25.1}
It follows from the proof that in the definition of $\phi_{t}$ one can take
 any constant
$\geq B$ in place of $B$.

\end{remark}

\mysection{Proof of Theorem \protect\ref{theorem 12.24.1}}
                       \label{section 1.14.1}

We prove the theorem after some discussion.

\begin{lemma}
                         \label{lemma 12.8.1}

For any $n\geq1$,  $p_{i} >0$, $i=1,...,n$,  and
$\kappa\geq d+1+p_{1}+...+p_{n}$ there exists a constant
$N=N(d,\kappa,n,p_{i})$ such that
for  any     $A_{i}\in \bR^{d}$, $i=1,...,n$, we have
\begin{equation}
                        \label{12.18.1}
|A_{1}|^{p_{1}}\cdot...\cdot|A_{n}|^{p_{n}}\leq 
N\int_{\bR^{d}}h(\eta)|(\eta,A_{1})|^{p_{1}}\cdot...\cdot |(\eta,A_{n})|^{p_{n}}
\,d\eta,
\end{equation}
where $h(\eta)=(1+|\eta|^{\kappa})^{-1}$
and $(\eta,A_{i})$ is the scalar product
of $\eta$ and $A_{i}$.
\end{lemma}

Proof. By dividing both parts by $|A_{1}|^{p_{1}}$
we see that we may assume that $|A_{1}|=1$.
After that we may also assume that $|A_{i}|=1$
for all $i$.
Then the result easily follows by contradiction.
The lemma is proved.

By applying this lemma to $A_{1}=Du$ and $A_{2}=Dv_{x^{i}}$ and then summing up with respect to $i$ we get the following, where
$$
u_{(\eta)}:=\eta^{i}D_{i}u.
$$

\begin{corollary}
                    \label{corollary 12.18.2}
For any smooth functions $u,v$ on $\bR^{d}$,
on $  \bR^{d}$ we have for appropriate $p,q,\kappa$
$$
|Du|^{p} 
\leq N\int_{\bR^{d}}h(\eta) |u_{(\eta)}|^{p} \,d\eta,
$$
$$
|Du|^{p}|D^{2}v|^{q}
\leq N\int_{\bR^{d}}h(\eta) |u_{(\eta)}|^{p}|Dv_{(\eta)}|^{q} \,d\eta.
$$

\end{corollary}

{\bf Proof of  Theorem \ref{theorem 12.24.1}}.
Fix $\kappa\geq d+1+p$ and for functions $u(x,\eta)$ and $v(x,\eta)$
on $\bR^{2d}$ let us write $u\prec_{\kappa} v$
if 
$$
\int_{\bR^{2d}}h(\eta)  u(x,\eta) \,dxd\eta
\leq \int_{\bR^{2d}}h(\eta)  v(x,\eta) \,dxd\eta.
$$
 We also write
$u\sim_{\kappa} v$ if the above  integrals coincide. In these terms Corollary \ref{corollary 12.18.2} implies that
\begin{equation}
                        \label{12.18.2}
|\eta|^{p } |Du|^{p} \prec_{\kappa} N|u_{(\eta)}|^{p} ,\quad
 |\eta|^{p} |Du|^{p-2}|D^{2}v|^{2}\prec_{\kappa} N|u_{(\eta)}|^{p-2}|Dv_{(\eta)}|^{2}.
\end{equation}
 
Introduce
$$
v_{t}=v_{t}(x,\eta)=\eta^{i}D_{i}
u_{t}(x),\quad \theta_{t}=\sign v_{t},\quad 
V_{t}=|v_{t}|^{p/2},
$$
$$
 W^{2}_{t}:=|\eta|^{p}|Du_{t}|^{p-2}
|D^{2}u_{t}|^{2},\quad\gamma^{j}_{t}=D_{i}a^{ij}_{t}+b^{i}_{t}.
$$
Observe that thanks to \eqref{12.18.2}
\begin{equation}
                           \label{12.18.5}
|\eta|^{p } |Du|^{p} \prec_{\kappa} N|v_{t}|^{p},\quad
 W^{2}_{t}\prec_{\kappa}N |D V_{t}|^{2},\quad |D V_{t}|^{2}\leq N  W^{2}_{t}.
\end{equation}

Owing to \eqref{1.10.6} we can substitute
$\phi_{(\eta)}$ in place of $\phi$
into \eqref{1.11.2} which yields
$$
dv_{t}
= (L_{t}v_{t}+a^{ij}_{t(\eta)}D_{ij}u_{t}+
\gamma^{i}_{t(\eta)} D_{i}u_{t}
+c_{t(\eta)}u_{t}+f_{t(\eta)})\,dt
$$
$$
+(M^{k}_{t}v_{t }+
\sigma^{ik}_{t(\eta)}
D_{i}u_{t}+\nu^{k}_{t(\eta)}u_{t}+g^{k}_{t(\eta)})\,dw^{k}_{t}.
$$

Again in light of Lemma \ref{lemma 1.6.1},
Theorem 2.1 of \cite{Kr_10} is applicable
to $v_{t}$,
which in particular implies the $L_{p}$-continuity of $v_{t\wedge\tau}$ and thus
of $Du_{t\wedge\tau}$. By this theorem
also  we obtain for $t\leq\tau$
  that
$$
(1/p) \int_{\bR^{2d}}h(\eta)|v_{t}|^{p}\,dxd\eta=
(1/p) \int_{\bR^{2d}}h(\eta)|v_{0}|^{p}\,dxd\eta
$$
$$
+\int_{0}^{t}\int_{\bR^{2d}}h(\eta)\theta_{s}|v_{s}|^{p-1}\big(L_{s}v_{s}+a^{ij}_{s(\eta)}D_{ij}u_{s}\big)\,dxd\eta\,ds
$$
$$
 +\int_{0}^{t}\int_{\bR^{2d}}h(\eta)\theta_{s}|v_{s}|^{p-1}\big(\gamma^{i}_{t(\eta)} D_{i}u_{t}
+c_{s(\eta)}u_{s}+f_{s(\eta)}\big)\,dxd\eta\,ds 
$$
$$
+ 
\int_{0}^{t}\int_{\bR^{2d}}h(\eta)J_{s}\,dxd\eta\,ds+ m_{t},
$$
where   
$$
J_{t}=((p-1)/2)|v_{t}|^{p-2}\sum_{k }
\big(M^{k }_{t}v_{t }+
\sigma^{ik }_{t(\eta)}
D_{i}u_{t}+\nu^{k }_{t(\eta)}u_{t}+g^{k }_{t(\eta)}\big)^{2}
$$
$$
\leq ((p-1)/2)|v_{t}|^{p-2}\sigma^{ik}_{t}D_{i}v_{t}\sigma^{jk}_{t}D_{j}v_{t}+(\delta/(2p))|DV_{t}|^{2}
$$

\begin{equation}
                                \label{12.22.1}
+N|\nu_{t}|^{2}|v_{t}|^{p} +N|\eta|^{2}|v_{t}|^{p-2}\big(|D\sigma_{t}|^{2}|Du_{t}|^{2}+ |D\nu_{t}|^{2}|u_{t}|^{2}+ |Dg_{t}|^{2}\big),
\end{equation}
$$
 m_{t}=\int_{0}^{t}\int_{\bR^{2d}}h(\eta)J^{k}_{s}\,dxd\eta\,dw^{k}_{s},
$$
 $$
J^{k}_{t}=\theta_{t}|v_{t}|^{p-1}
\big(M^{k }_{t}v_{t }+
\sigma^{ik }_{t(\eta)}
D_{i}u_{t}+\nu^{k }_{t(\eta)}u_{t}+g^{k }_{t(\eta)}\big).
$$

Next, we may assume that $\widehat{Da},\widehat{D\sigma},\hat b,\hat\cedc,
\hat\frc,\hat\nu,\widehat{D\nu}\leq 1$ and
first   deal with the terms  
containing
$a$ and $\sigma$ but not their derivatives.
Their sum is (see \eqref{12.22.1}) less than
$$
I_{1}:=\theta_{t}|v_{t}|^{p-1}D_{i}(a^{ij}_{t}D_{j}v_{t})
+((p-1)/2)|v_{t}|^{p-2}\sum_{k }(\sigma^{ik }D_{i}v_{t})^{2} 
$$
$$
+(\delta/(2p))|D V_{t}|^{2}
\sim_{\kappa}-((p-1) /p^{2})
(4a^{ij}_{t}-2\alpha^{ij}_{t})D_{i} V_{t}D_{j} V_{t}+(\delta/(2p))|D V_{t}|^{2}
$$
\begin{equation}
                               \label{12.20.2}
\leq - (\delta/(2p))|D V_{t}|^{2}\prec_{\kappa} -
\hat \delta(|D V_{t}|^{2}+ W^{2}_{t}),
\end{equation}
where $\hat \delta=\hat \delta(p,d,\delta)>0$
and in the last inequality we used \eqref{12.18.5}.

By Lemma \ref{lemma 12.17.1} 
$$
I_{2}:=\theta_{t}v_{t}^{p-1}a^{ij}_{t(\eta)}D_{ij}u_{t}
\leq ( V_{t}|Da_{t}|)|\eta|^{p/2}|Du_{t}|^{p/2-1}
|D^{2}u_{t}|
$$
$$
\leq(\hat\delta/2) W^{2}_{t}+N  V_{t}^{2}|Da_{t}|^{2}\prec_{\kappa}(\hat\delta/2)  W^{2}_{t} 
+N_{1}\widehat{Da}^{2}|D V_{t}|^{2}+N(1+\overline{Da}_{t }^{2})
 V_{t}^{2}.
$$
We subject $\widehat{Da}$ to
\begin{equation}
                           \label{12.24.3}
N_{1}\widehat{Da}^{2}\leq \hat\delta/2
\end{equation}
and get that

\begin{equation}
                           \label{12.18.4}
I_{2} \prec_{\kappa} (\hat \delta/2) (|D V_{t}|^{2}+ W^{2}_{t})+N(1+\overline{Da}_{t}^{2}) V^{2}_{t}.
\end{equation}

Next, come  the  terms containing the products of $\gamma,b$ and their derivatives:
$$
I_{3}:=\theta_{t}|v_{t}|^{p-1}b^{i}_{t}D_{i}v_{t}+
\theta_{t}|v_{t}|^{p-1}\gamma^{i}_{t(\eta)} D_{i}u_{t} \sim_{\kappa}
(2/p)\big( V_{t}b^{i}_{t}\big)D_{i} V_{t}
$$
$$
-\eta^{j}\gamma^{i}_{t} \big((p-1)|v_{t}|^{p-2}D_{i}v_{t}
D_{i}u_{t}+|v_{t}|^{p-1}D_{ij}u_{t}\big).
$$ 
Here by Corollary \ref{corollary 12.17.1}  
$$
(2/p)\big( V_{t}b^{i}_{t}\big)D_{i} V_{t}\prec_{\kappa}  (\hat\delta/8)|D V_{t}|^{2}+
N_{5}\hat b^{2}|D V_{t}|^{2}+N\bar b_{t}^{2} V_{t}^{2}.
$$
By Corollary \ref{corollary 12.18.2} the remaining term in $I_{3}$ is 
dominated  by   
$$
N  |\gamma_{t}| |\eta|^{p}|Du|^{p-1}|D^{2}u|
$$
$$\prec_{\kappa}
N(|Da_{t}|+|b_{t}|) V_{t}\big(|\eta|^{p/2}|Du_{t}|^{p/2-1}|D^{2}u_{t}|\big)
$$
$$
\leq N(|Da_{t}|^{2}+|b_{t}|^{2}) V_{t}^{2}+(\hat \delta/4) W^{2}_{t}
$$
$$
\prec_{\kappa} N_{2}(\widehat{Da}^{2}+\hat b^{2})|D V_{t}|^{2}+
(\hat \delta/4) W^{2}_{t}+N(1+\overline{Da}^{2}_{t}+\bar b^{2}_{t})
 V_{t}^{2}.
$$
 
We subject $\widehat{Da}$ and  $\hat b$ to 
\begin{equation}
                           \label{12.24.4}
N_{2}(\widehat{Da}^{2}+\hat b^{2})\leq \hat\delta/8,
\end{equation}
and  get that
\begin{equation}
                           \label{12.19.02}
I_{3} \prec_{\kappa} (\hat \delta/4) (|D V_{t}|^{2}+ W^{2}_{t})+N(1+\overline{Da}^{2}_{t}+\bar b_{t}^{2}) V^{2}_{t}.
\end{equation}

The term with $|Du_{t}|^{2}$ is 
$$
I_{4}:=N|\eta|^{2}|v_{t}|^{p-2}|D\sigma_{t}|^{2}|Du_{t}|^{2}\leq N|\eta|^{p}|D\sigma_{t}|^{2}
\big(|Du_{t}|^{p/2}\big)^{2}
$$
$$
\prec_{\kappa}
N_{3} \widehat{D\sigma} ^{2} W^{2}_{t}+N(1+\overline{D\sigma}_{t}^{2}) V^{2}_{t}.
$$
We subject $\widehat{D\sigma}$ to 
\begin{equation}
                           \label{12.24.5}
N_{3} \widehat{D\sigma}^{2}\leq \hat\delta/8,
\end{equation}
and  get that
\begin{equation}
                           \label{12.18.9}
I_{4} \prec_{\kappa} (\hat \delta/8) (|D V_{t}|^{2}+ W^{2}_{t})+N(1+\overline{D\sigma}_{t}^{2}) V^{2}_{t}.
\end{equation}

Now come the terms which do not contain $Du_{t}$.
Their sum is
$$
I_{5}:=I_{6}+I_{7}+I_{8},
$$
where
$$
I_{6}:=\theta_{t}|v_{t}|^{p-1}\big(c_{t}v_{t}+c_{t(\eta)}u_{t}\big),
\quad
I_{7}:=N|\nu_{t}|^{2}|v_{t}|^{p}
+N|\eta|^{2}|v_{t}|^{p-2}|D\nu_{t}|^{2}
|u_{t}|^{2},
$$
$$
I_{8}:=\theta_{t}|v_{t}|^{p-1}f_{t(\eta)}+
N|\eta|^{2}|v_{t}|^{p-2}|Dg_{t}|^{2}.
$$

We have
$$
I_{6} \leq \cedc_{t}^{2} V_{t}^{2}+\frc_{t}^{2}
\big( V_{t}^{2}+|\eta|^{p}|u_{t}|^{p}\big)\prec_{\kappa}
N_{4}\big( \hat \cedc^{2}+\hat\frc^{2}\big)|D V_{t}|^{2} 
$$
$$
+N\hat\frc^{2}|\eta|^{p} |DU_{t} |^{2}+N(
(1+\bar \frc_{t}^{2})|\eta|^{p}|u_{t}|^{p}+(1+\bar\cedc_{t}^{2}+\bar \frc_{t}^{2})
 V_{t}^{2})
$$

We subject $\hat \cedc^{2}$ and $\hat\frc^{2}$ to
\begin{equation}
                           \label{12.24.6}
N_{4}\big(\hat \cedc^{2}+\hat\frc^{2}\big)\leq
\hat \delta/16
\end{equation}
and get that
$$
I_{6}\prec_{\kappa} (\hat \delta/16) (|D V_{t}|^{2}+ W^{2}_{t})
$$
\begin{equation}
                           \label{12.24.1}
+N(1+\bar\cedc_{t}^{2}+\bar \frc_{t}^{2})
 V_{t}^{2}+N |\eta|^{p} |DU_{t} |^{2} +N
(1+\bar \frc_{t}^{2})|\eta|^{p}U_{t}^{2}.
\end{equation}

Similarly we deal with $I_{7}$ and see that,
if for some well defined $N_{5}$
\begin{equation}
                           \label{12.24.7}
N_{5}(\hat \nu^{2}+\widehat{D\nu}^{2})\leq
\hat\delta/32,
\end{equation}
then
$$
I_{7}\prec_{\kappa} (\hat \delta/32) (|D V_{t}|^{2}+ W^{2}_{t})
$$
\begin{equation}
                           \label{12.24.2}
+N(1+\bar\nu_{t}^{2}+\overline{D\nu}_{t}^{2})
 V_{t}^{2}+N |\eta|^{p} |DU_{t} |^{2} +N
(1+\overline{D\nu}_{t}^{2})|\eta|^{p}U_{t}^{2}.
\end{equation}

In what concerns $I_{8}$, obviously,
$$
I_{8}\leq V_{t}^{2}+N|\eta|^{p}(|Df_{t}|^{p}
+|Dg_{t}|^{p}).
$$

Thus, given that \eqref{12.24.3}, \eqref{12.24.4}, 
\eqref{12.24.5}, \eqref{12.24.6}, and
\eqref{12.24.7} hold, we have  
$$
d|v_{t}|^{p}\prec_{\kappa}
-\check\delta (|D V_{t}|^{2}+ W^{2}_{t})
\,dt+(C-1)\Lambda_{t}|v_{t}|^{p}\,dt
$$
\begin{equation}
                          \label{12.24.8}
+N|\eta|^{p}\big(|Df_{t}|^{p}
+|Dg_{t}|^{p}+  |DU_{t} |^{2} + 
(1+\bar\frc_{t}^{2}+\overline{D\nu}_{t}^{2}) U_{t}^{2}\big)\,dt
+J^{k}_{t}\,dw^{k}_{t},
\end{equation}
where $\check\delta=\hat \delta/(32p)$.

It follows that
$$
d\Big(e^{-\Psi_{t}}|v_{t}|^{p}\Big)\prec_{\kappa}
-\check\delta e^{-\Psi_{t}}(|D V_{t}|^{2}+ W^{2}_{t})
\,dt -\Lambda_{t}e^{-\Psi_{t}}|v_{t}|^{p}\,dt
$$
\begin{equation}
                               \label{1.13.2}
+Ne^{-\Psi_{t}}|\eta|^{p}\big(|Df_{t}|^{p}
+|Dg_{t}|^{p}+  |DU_{t} |^{2} + 
(1+\bar\frc_{t}^{2}+\overline{D\nu}_{t}^{2}) U_{t}^{2}\big)\,dt
+e^{-\Psi_{t}}J^{k}_{t}\,dw^{k}_{t}.
\end{equation}
We convert  this into the
 integral form with integrals with respect to
$x,\eta$, and $t$, then
take expectations of the resulting inequality.
We also use Corollary \ref{corollary 12.18.2}.
Then we get
$$
H=E\int_{0}^{\tau}e^{-\Psi_{t}}\int_{\bR^{d}}\big(
|Du_{t}|^{p-2}
|D^{2}u_{t}|^{2}+\Lambda_{t}|Du_{t}|^{p}\big)\,dxdt
\leq NE|Du_{0}|^{p}
$$
\begin{equation}
                              \label{1.13.3}
+
NE\int_{0}^{\tau}e^{-\Psi_{t}}
\int_{\bR^{d}} \big(|Df_{t}|^{p}
+|Dg_{t}|^{p}+  |DU_{t} |^{2} + 
(1+\bar\frc_{t}^{2}+\overline{D\nu}_{t}^{2}) U_{t}^{2}\big)\,dt.
\end{equation}
We may assume that $C 
\Lambda_{t}\geq B\lambda_{t}+
1+\bar\frc_{t}^{2}+\overline{D\nu}_{t}^{2}$
and then, owing to \eqref{12.19.2},
  the term in \eqref{1.13.3} containing 
$|DU_{t} |^{2} + 
(1+\bar\frc_{t}^{2}+\overline{D\nu}_{t}^{2}) U_{t}^{2}$ is estimated by
$$
NE|u_{0}|^{p}+NE\int_{0}^{\tau}e^{-\Psi_{t}}
\int_{\bR^{d}}\big(|f_{t}|^{p}+|g_{t}|^{p}\big)\,dxdt\leq NK,
$$
so that $H\leq NK$.

Next, we come back to \eqref{1.13.2}
and, after converting it,  
take sup
with respect to $t\leq\tau$ and take expectations in the resulting inequality.
We also use Corollary \ref{corollary 12.18.2}
and the result of the above treatment of
\eqref{1.13.3}.
Then we obtain
\begin{equation}
                          \label{12.24.80}
G\leq NK+I  ,
\end{equation}
where
$$
I=E\sup_{t\leq\tau}\Big|\int_{0}^{t}e^{-\Psi_{s}}
\int_{\bR^{2d}}h(\eta)J^{k}_{s}\,dxd\eta\,dw^{k}_{s}\Big|.
$$
 
By Davis's inequality
$$
I  
\leq NE\Big(\int_{0}^{\tau}e^{-2\Psi_{t}}\sum_{k}\Big(\int_{\bR^{2d}}h(\eta)J^{k}_{t}\,dxd\eta\Big)^{2}\,dt\Big)^{1/2}
$$
$$
= NE\Big(\int_{0}^{\tau}\sum_{k}e^{-2\Psi_{t}}\Big(
\int_{\bR^{2d}}h(\eta)\theta_{t}|v_{t}|^{p-1}
\big(M^{k }_{t}v_{t }
$$
$$
+
\sigma^{ik }_{t(\eta)}
D_{i}u_{t}
+\nu^{k }_{t(\eta)}u_{t}+g^{k }_{t(\eta)}\big)\,dxd\eta\Big)^{2}\,dt\Big)^{1/2}
$$
$$
\leq NE\Big(\sup_{t\leq\tau}e^{-\Psi_{t}}
\int_{\bR^{d}}|v_{t}|^{p}\,dx\Big)^{1/2}\Big(\int_{0}^{\tau}e^{-\Psi_{t}}\int_{\bR^{2d}}h(\eta)J_{t}\,dxd\eta dt\Big)^{1/2}
$$
$$
\leq NG^{1/2}\Big(E\int_{0}^{\tau}e^{-\Psi_{t}}\int_{\bR^{2d}}h(\eta)J_{t}\,dxd\eta dt\Big)^{1/2}.
$$

Above we actually estimated $J_{t}$ splitting it in parts. It is seen
from these estimates that
$$
J_{t}\prec_{\kappa} N|DV_{t}|^{2}+I_{4}+I_{7}
+N|\eta|^{2}|v_{t}|^{p-2}|Dg_{t}|^{2}
$$
$$
\prec_{\kappa}
 N|\eta|^{p}\big(|Du_{t}|^{p-2}
|D^{2}u_{t}|^{2}+\Lambda_{t}|Du_{t}|^{p}\big)
+ N |\eta|^{p} |DU_{t} |^{2} 
$$
$$
+N
(1+\overline{D\nu}_{t}^{2})|\eta|^{p}U_{t}^{2}
+N|\eta|^{p}|Dg_{t}|^{p}.
$$
Recalling our argument about \eqref{1.13.3}
we see that
 $
I\leq NG^{1/2}K^{1/2}$, which after being
substituted into \eqref{12.24.80},
yields that $G\leq NK$ as well.
 This proves \eqref{12.24.9} and the theorem.

The manuscript contains no data and there is no
conflict of interest.

\end{document}